\documentclass[12pt]{amsart}

\usepackage{times}
\usepackage[latin1]{inputenc}
\usepackage[T1]{fontenc}
\usepackage{amsmath,amssymb,amsthm}
\usepackage{graphicx}
\usepackage{epsfig}
\usepackage{hyperref}

\usepackage[square,authoryear,sort&compress]{natbib}

\begin{document}

\newtheorem{theorem}{Theorem}[section]
\newtheorem{lemma}[theorem]{Lemma}
\newtheorem{conjecture}[theorem]{Conjecture}
\newtheorem{question}[theorem]{Question}
\newtheorem{corollary}[theorem]{Corollary}
\newtheorem{proposition}[theorem]{Proposition}
\theoremstyle{definition}
\newtheorem{definition}[theorem]{Definition}
\newtheorem{example}[theorem]{Example}
\newtheorem{claim}[theorem]{Claim}
\newtheorem{xca}[theorem]{Exercise}

\newcommand{\ind}{\mbox{ind}}
\newcommand{\conn}{\mbox{connectivity}}
\newcommand{\Hom}{\mbox{Hom}}

\theoremstyle{remark}
\newtheorem{remark}[theorem]{Remark}
\newcommand{\be}{\begin{equation}}
\newcommand{\ee}{\end{equation}}
\newcommand{\lexmin}{\mbox{lexmin}}
\newcommand{\prob}{\mbox{\bf P}}
\newcommand{\nerve}{\mathcal{N}}
\newcommand{\Star}{\mbox{st}}
\newcommand{\supp}{\mbox{supp}}
\newcommand{\vsupp}{\mbox{vsupp}}
\newcommand{\link}{\mbox{lk}}
\newcommand{\Z}{\mathbb{Z}}
\newcommand{\R}{\mathbb{R}}
\newcommand{\Q}{\mathbb{Q}}
\newcommand{\Homology}{\widetilde{H}}
\newcommand{\bound}{\partial}
\newcommand{\XG}{X(n,p)}
%\theoremstyle{remark}
%\newtheorem{remark}[theorem]{Remark}

%    Blank box placeholder for figures (to avoid requiring any
%    particular graphics capabilities for printing this document).
%\newcommand{\blankbox}[2]{%
%  \parbox{\columnwidth}{\centering
%%    Set fboxsep to 0 so that the actual size of the box will match the
%%    given measurements more closely.
%    \setlength{\fboxsep}{0pt}%
%    \fbox{\raisebox{0pt}[#2]{\hspace{#1}}}%
%  }%
%}

\bibliographystyle{plain}

\title{Points in a triangle forcing small triangles \\}

\author{Matthew Kahle}
\address{Department of Mathematics, Stanford University}
\email{mkahle@math.stanford.edu}
\thanks{Supported in part by Stanford's NSF-RTG grant}

\date{November 14, 2008}

\maketitle

\vspace*{-0.5in}

\begin{abstract} An old theorem of Alexander Soifer's is the following: Given five points in a triangle of unit area, there must exist some three of them which form a triangle of area $1/4$ or less. It is easy to check that this is not true if ``five'' is replaced by ``four'', but can the theorem be improved in any other way? We discuss in this article two different extensions of the original result.

First, we allow the value of ``small'' $1/4$ to vary. In particular, our main result is to show that given five points in a triangle of unit area, then there must exist some three of them determining a triangle of area $6/25$ or less. %We also give an example of five points in a unit triangle where all the triangles formed are of area $\ge 1/6$. %, and we can not rule out the possibility that $A=1/6$.

Second, we put bounds on the minimum number of small triangles determined by $n$ points in a triangle, and make a conjecture about the asymptotic right answer as $n \to \infty$.

\end{abstract}

\section{Introduction} This article is in honor of Alexander ``Sasha'' Soifer's sixtieth birthday. I have known Sasha for about twenty years now; I first met him when I competed in his Colorado Math Olympiad in 1988. I remember that he teased me about looking even younger than I was (and also about my poor grades in school), but when we talked about mathematics he always spoke to me as his peer. We could speak a common language and share our curiosity and passion for solving problems.

%I took a class in combinatorial geometry from Alex, while I was still in high school, and I still remember it as one of my favorite classes. Many years later, my research tastes lean heavily toward combinatorics and geometry. I sometimes have wondered if this was due to Alex's early influence, and I am sure to some extent it is. But I know that when I was nine years old, I was completely obsessed with the Rubik's Cube, and what is it if not combinatorial and geometric.

 Professor Soifer has solved many outstanding math problems, but posed many more. He has composed most of the problems for twenty-five (and counting) Colorado Math Olympiads, and written numerous research articles and five books about open problems in combinatorics and geometry. Perhaps the most important thing that I learned from him is that we are free to ask our own mathematical questions and pursue them as far as we are able. This was a very empowering idea to me when I was young. %Many years later, as a working mathematician, I appreciate even more the importance of asking the right questions.

 Sasha has nice taste in problems too, and tells a good mathematical story \cite{Coloring}; I have sometimes felt a little haunted by problems which I first heard from him. One problem I have always remembered for its simple and elegant statement (and somewhat tricky proofs) is the following, which appeared on the Colorado Math Olympiad in 1988 \cite{Olympiad}. It stayed with me, and I have occasionally sat down trying to see a new solution to it.

\begin{theorem} \label{soif} Given five points in a triangle of unit area, there must exist some three of them forming a triangle of area $1/4$ or less. \end{theorem}

If the reader has never seen this problem, they are warmly invited to try it. The problem becomes easier if ``five'' is replaced by ``six'' or ``seven,'' but it seems to take some work to get ``five''. Solutions by Alexander Soifer, Royce Peng, and Cecil Rosseau all appear in the book, ``How does one cut a triangle?'' which will soon be reprinted by Springer \cite{HowDoes}. One wonders if this is the best possible result, and in a sense it is: it is easy to check that ``five'' can not be replaced by ``four''. But what about the $1/4$?

I must admit that I did not make any real progress toward proving Theorem \ref{soif} during the Olympiad, but now I will have a (very) small revenge on the problem by improving $1/4=0.25$ to $6/25=0.24$. In the spirit of Sasha Soifer and his wonderful journal Geombinatorics, we leave the reader with some open problems as well.

\section{Varying the definition of small}

Our main result is the following. %{ \bf We use the word ``small'' to refer to triangles of area $\le 6/25=24/100$ for the rest of the section.} %We rely on ideas from proofs of Theorem \ref{soif} due to Peng and Soifer.

\begin{theorem} Given five points in a triangle of unit area, then some three of them form a triangle of area less than or equal to $6/25$.
\end{theorem}

\begin{figure}
 \includegraphics[width=2.25in]{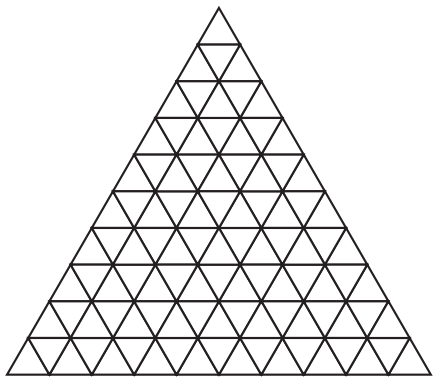}
\caption{Partitioning $T$ into $100$ congruent triangles.}
\label{fig:hund}
\end{figure}

Since ratios of areas are preserved under affine transformation, we assume without loss of generality that our unit triangle $T$ is equilateral. Let $L$ denote the length of one side of $T$ (i.e. $L=2/\sqrt[4]{3}$). We partition $T$ into $100$ congruent equilateral triangles of side length $L/10$, as in Figure \ref{fig:hund}. Label the five points by $P=\{p_1, p_2, \ldots, p_5\}$.

{ \bf We use the word ``small'' to refer to triangles of area $\le 6/25$ for the rest of the section.}

\begin{lemma}\label{dist4} If for any pair of points in $p_i, p_j \in P$, $d(p_i, p_j) \le 4L/10$ then three of the points in $S$ form a small triangle. (In particular, if $d(p_1,p_2) < 4L / 10$ then either $\{ p_1, p_2, p_k \}$ forms a small triangle for some $k \in  \{3, 4,5\}$, or else $\{p_3, p_4, p_5\}$ forms a small triangle.)
\end{lemma}

\begin{proof}[Proof of Lemma \ref{dist4}] For any three points $x,y,z$ let $A(x,y,z)$ denote the area of the triangle they form. The locus of points $x$ such that $A(p_i,p_j,x) \le 24/100$ is a closed infinite strip $S$ centered on line $p_ip_j$, and of width at least $$2 ( \sqrt{3} /2)(6L/10)=(3\sqrt{3}/5)L > L.$$ Since the diameter of $T$ is $L$, the intersection $I$ of the complement of the strip with $T$, has at most one component. The main point is that if a third point is inside $S$ we are done, and otherwise all three points lie in $I$, and we need only check that this forces the area of the triangle they form to be less than $24/100$.

%If $I$ is empty we are done, and if it is nonempty then its boundary is a polygon, either a triangle or a quadrilateral. The maximum area of a triangle with all three of its vertices in a closed polygon $P$ is realized by placing all vertices of the triangle at vertices of $P$.

By moving $p_i$ and $p_j$ along line $p_ip_j$ until they intersect the boundary of $T$ if necessary, we only increase the width of the strip $S$, so we may assume without loss of generality that $p_i$ and $p_j$ are both on the boundary of $T$.
If $p_i$ and $p_j$ are on the same side of $T$ then $I$ is a triangle with maximal area when $d(p_i,p_j)=4L/10$, and this area is only $(4/10)^2=16/100$, in which case we are done, so we may assume that $p_i$ and $p_j$ are on two different sides of $T$.

Now suppose that $d(p_i,p_j) < 4L / 10$. Then we can replace $p_i$ and $p_j$  by $p_i'$ and $p_j'$, respectively, so that the line $p_i'p_j'$ is parallel to $p_ip_j$, and $d(p_i', p_j') = 4L / 10$. This results in a new strip $S'$ and denote the intersection of $S'$ with $T$ by $I'$. We claim that $I'$ strictly contains $I$, as follows. Suppose that the distance between line $p_i p_j$ and $p_i' p_j'$ is $\Delta$. The strip $S$ is parallel to the strip $S'$, and the width $w'$ of strip $S'$ is related to the width $w$ of $S$ by $$(4L/10)w'=d(p_i,p_j) w.$$ We also have by basic trigonometry that $$4L/10 - d(p_i,p_j) = \Delta \left( \frac{1}{\tan{A}} + \frac{1}{\tan{2 \pi /3 - A}} \right),$$ where $A$ is the angle that line $p_ip_j$ makes with the sides of $T$, at $p_j$, as in Figure \ref{fig:seglemma}.

\begin{figure}
 \includegraphics[width=4in]{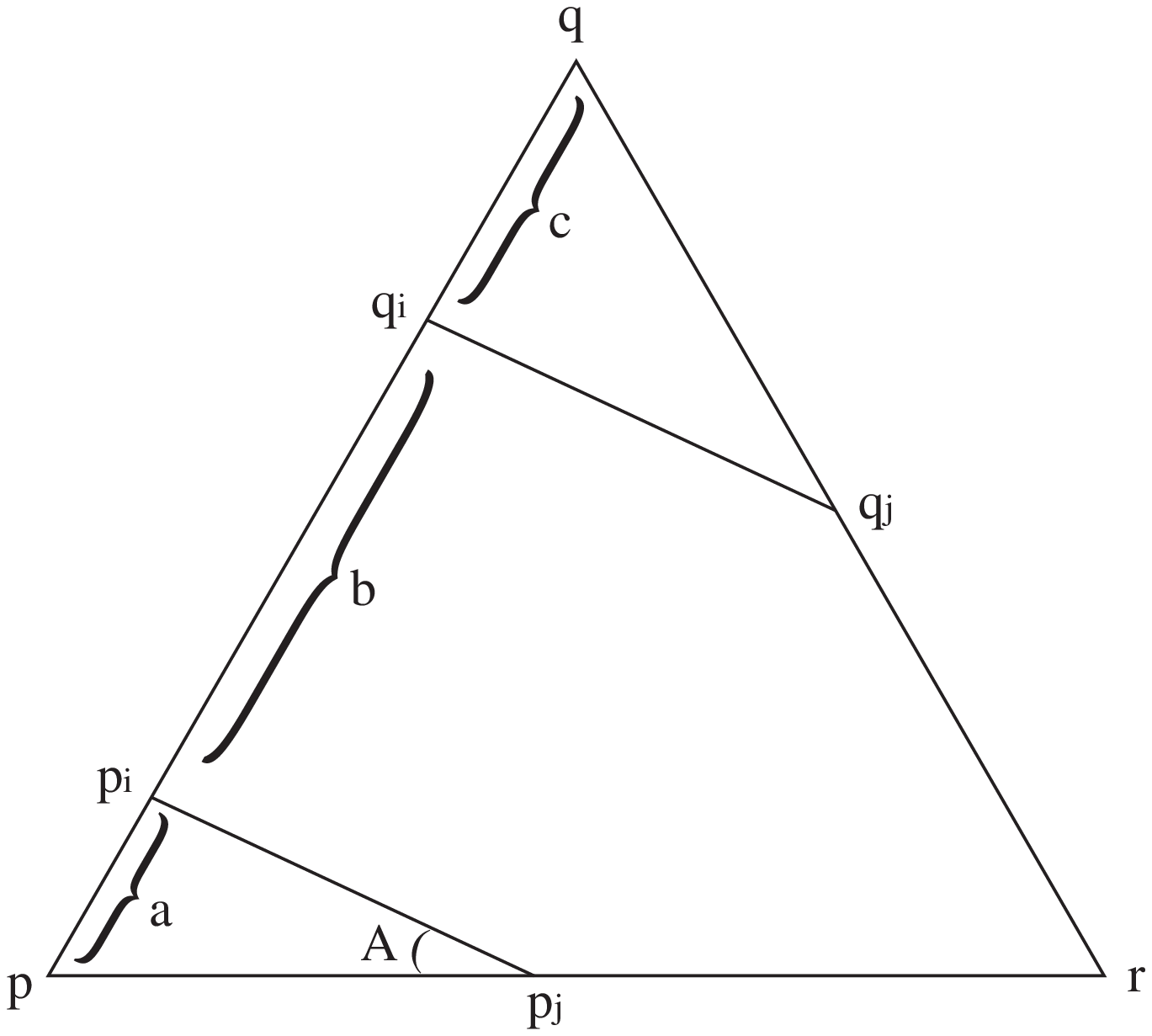}
\caption{Proof of Lemma \ref{dist4}.}
\label{fig:seglemma}
\end{figure}

The minimum of the function $f(x) = 1 / \tan(x) + 1 / \tan(2\pi/3 - x)$ on the interval $x \in [0, 2 \pi /3 ]$ is attained when $x= \pi/3$ and $f(x)= 2/ \sqrt{3} > 1$. Now

\begin{eqnarray*}
w-w' & = & w \left( 1- \frac{d(p_i,p_j) }{4L/10} \right) \\
&= & w \left( \frac{4L/10 - d(p_i,p_j)}{4L/10} \right) \\
& \ge & \frac{w}{4L/10}(2/\sqrt{3})\Delta\\
& \ge & \Delta,\\
\end{eqnarray*}
since $w \ge (3\sqrt{3} /10) L > 4L/10$.

This implies that $I'$ contains $I$, and we may assume that $d(p_i,p_j) = 4L /10$. 

Now label the corners of $T$ by $p$, $q$, and $r$, so that $p_i$ lies on side $pr$ and $p_j$ on side $pq$, as in Figure \ref{fig:seglemma}. If $d(p,p_i)$ and $d(p,p_j)$ are both greater than $4L/10$ then $d(p_i,p_j) > 4L/10$, and if they are both less than $4L/10$ then $d(p_i,p_j) < 4L/10$, but we are assuming that $d(p_i,p_j)=4L/10$. So by symmetry, we assume without loss of generality that $d(p,p_i) \le 4L / 10$ and $d(p,p_j) \ge 4L / 10$. In particular that $0 \le A \le \pi /3$.

Let $q_i$ and $q_j$ denote the other two intersections of the boundary of $S$ with the boundary of $T$. Note that line $q_iq_j$ intersects side $qr$, and in particular $I$ is a triangle. %Denote the measure of angle $p p_j p_i$ by $A$.
Let $a=d(p,p_i)$. Then by the Law of Sines, $$a= \sin{A(4L/10)}{\sqrt{3}/2}  =\frac{4 \sin{A}}{5\sqrt{3}} L.$$ Denote $b=d(p_i,q_i)$ and $c=d(q_i,r)$, so in particular $a+b+c=L$. We compute $b$ in terms of $A$ by
$$b=\frac{3\sqrt{3}/10}{\cos{(A-\pi/6)}},$$
so $$c= L \left(   1-\frac{4}{5\sqrt{3}} \sin{A} - \frac{3\sqrt{3}}{10\cos{(A-\pi/6)}}  \right).$$

The angles of triangle $q_iq_jq$ are $\pi/3$, $\pi /3 -A$, and $\pi/3 +A$, and applying the Law of Sines and a familiar formula for area of a triangle, its area is given by
\begin{eqnarray*}
|I| &= &\frac{1}{2} c^2 \frac{\sin{(\pi/3)}\sin{(\pi/3+A)}}{\sin{(\pi/3-A)}}.\\
\end{eqnarray*}

It seems that this is probably a monotone decreasing function of $A$ on the interval $A \in [0,2\pi/3]$ , but we did not care to take derivatives of this function to prove it. So we instead take a slightly more indirect approach, and one that does not involve calculus. The idea is to start with the fact that $c$ cannot be too large, which forces $A$ large, which in turn forces $c$ even smaller, etc., and we quickly reach a contradiction to the supposition that $|I| \ge 24/100$.

Since $0 \le A \le 2\pi /3$, we have $c \le 4L/10$, so $c^2 \le 16L^2/100$. Then if $|I| \ge 24/100$,
$$\frac{\sin{(\pi/3)}\sin{(\pi/3+A)}}{2\sin{(\pi/3-A)}} \ge  \frac{3}{2L^2},$$
so
$$\frac{\sin{(\pi/3+A)}}{\sin{(\pi/3-A)}} \ge  \frac{3}{2}.$$

Now $\sin{(\pi/3 + A )} \le 1$, so this implies in particular that $\sin{(\pi/3 -A)} < 2/3$, and then since $0 \le A \le \pi /3$ this gives that $0.317 < A$. Plugging back in to the formula for $c$ gives
\begin{eqnarray*}
c & = & L \left(   1-\frac{4}{5\sqrt{3}} \sin{A} - \frac{3\sqrt{3}}{10\cos{(A-\pi/6)}}  \right) \\
 & \le & L \left(   1-\frac{4}{5\sqrt{3}} \sin{A} - \frac{3\sqrt{3}}{10}  \right) \\
  & \le & L \left(   1-\frac{4}{5\sqrt{3}} \sin{0.317} - \frac{3\sqrt{3}}{10}  \right) \\
  & < & 0.337 L \\
 \end{eqnarray*}
  
We repeat the argument from before. If $c < 0.337L$ and $|I| \ge 24/100$, then $\sin{(\pi/3 -A)} \le c^2/0.24 < 0.474$, so $A > 0.553$. This gives in turn that 
\begin{eqnarray*}
c  & \le & L \left(   1-\frac{4}{5\sqrt{3}} \sin{0.553} - \frac{3\sqrt{3}}{10}  \right) \\
 & < & (24/100)L.\\
\end{eqnarray*}
a contradiction to the assumption that $|I| > 24/ 100$, since $$|I| = c d(r,q_j) / L^2 < c / L.$$

\end{proof}

Then we immediately have the following.

\begin{lemma}\label{trihex} If any two points of $P$ lie in a triangle of side length $4L/10$, or a hexagon of diameter $4L/10$, then some three points form a small triangle.
\end{lemma}

\begin{figure}
 \includegraphics[width=2in]{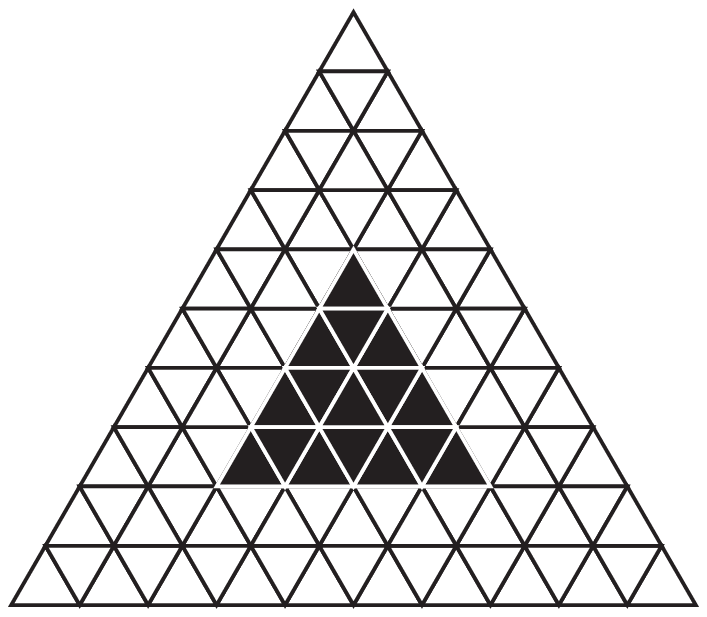}
 \includegraphics[width=2in]{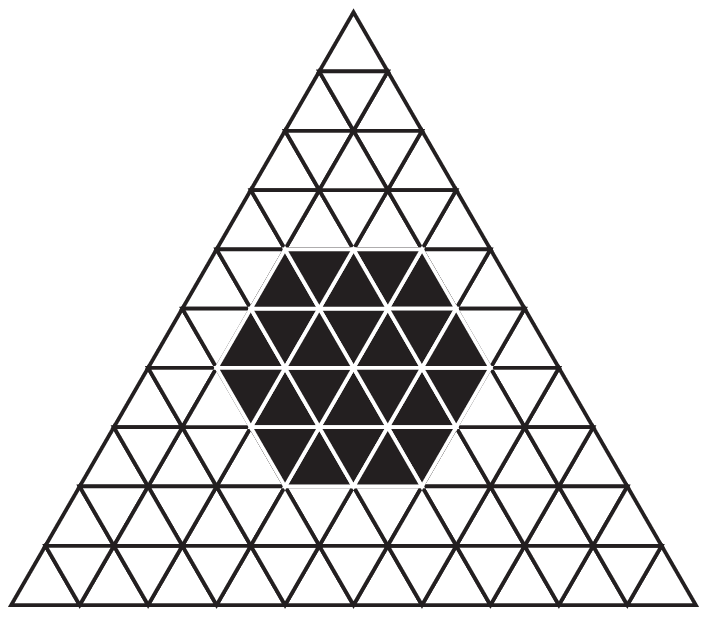}
\caption{A triangle and hexagon of diameter $4L/10$.}
\label{fig:trihex}
\end{figure}

These are many such triangles and hexagons suggested by our partition of $T$, as in Figure \ref{fig:trihex}.

We further restrict the possible arrangements of points by using the fact that a triangle in a parallelogram of area $A$ has area at most $A/2$. (This observation is used in Soifer's original proof of Theorem \ref{soif}.)

\begin{lemma} \label{parallel} If any three points of $P$ lie in a parallelogram of area $A \le 48/100$, then they form a triangle of area no more than $A/2 \le 24/100$.
\end{lemma}

\begin{figure}
 \includegraphics[width=1.5in]{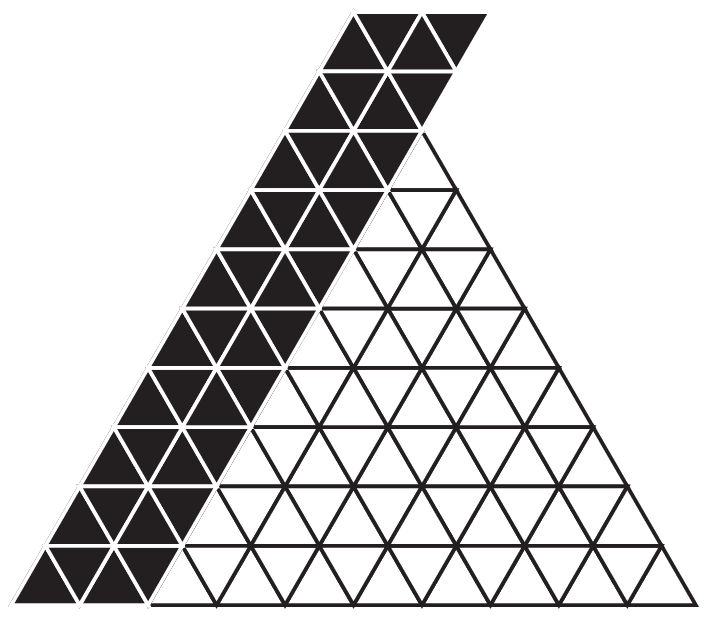}
 \includegraphics[width=1.5in]{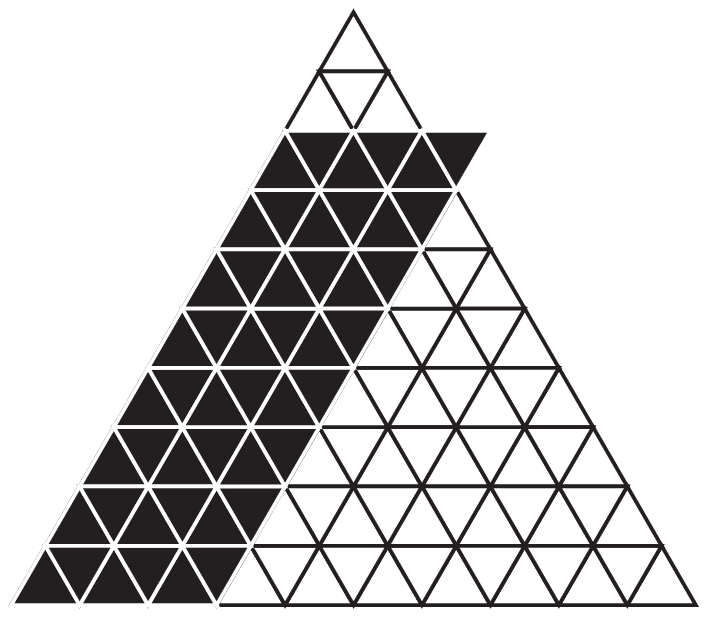}
  \includegraphics[width=1.5in]{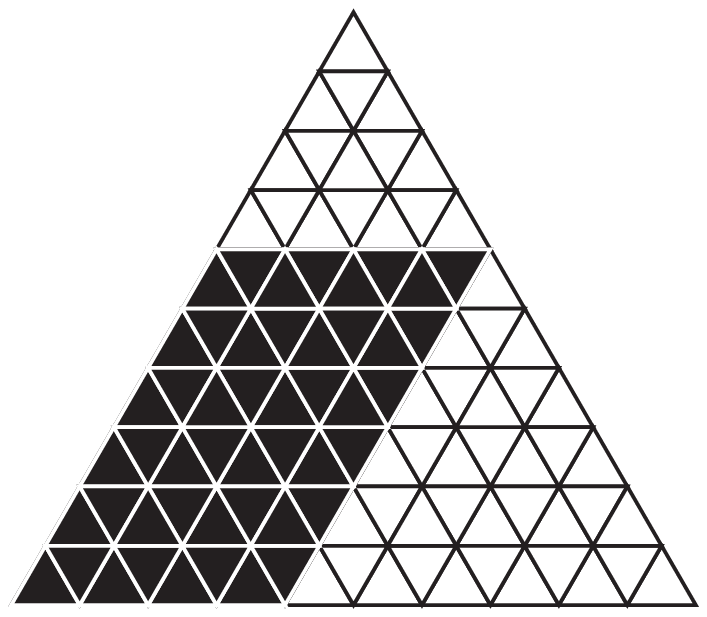}
 \caption{Parallelograms of area $A\le48/100$. From left to right, $(2,10)-$, $(3,8)-$, and $(4,6)-$ parallelograms.}
\label{fig:parallel}
\end{figure}

We make use of three such kinds of parallelograms, illustrated in Figure \ref{fig:parallel}, which we will call $(2,10)-$, $(3,8)-$, and $(4,6)-$ parallelograms.

We consider three cases.

\begin{figure}
\includegraphics[width=2in]{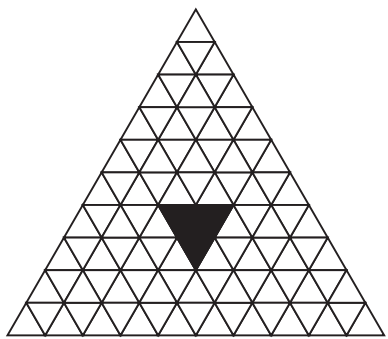}
\caption{Case I: At least one point is in the shaded triangle.}
\label{fig:case1}
\end{figure}

Case I is that there is a point of $P$ inside the central shaded triangle in Figure \ref{fig:case1}.

\begin{figure}
\includegraphics[width=2in]{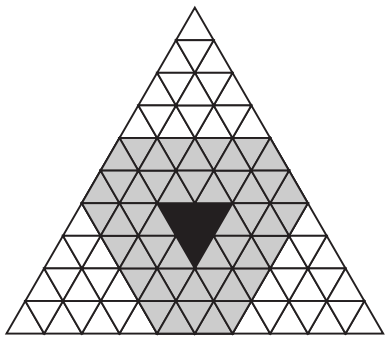}
\caption{Case I: There can be no more than the given point in the greater shaded region. Hence four points lie in the three white triangles.}
\label{fig:wide1}
\end{figure}

Assuming Case I, by the hexagon case of Lemma \ref{trihex}, there can be no other points in the wider shaded region in Figure \ref{fig:wide1}, or else we are already done. But that leaves at least four points in the three corner white triangles. Two of them must lie in the same triangle, but then by the triangle case of Lemma \ref{trihex} we are done.

\begin{figure}
\includegraphics[width=2in]{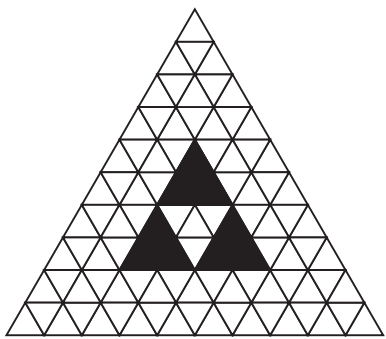}
\caption{Case II: At least one point $p_1$ lies in one of the three shaded triangles. Assume that $p_1$ in the top triangle.}
\label{fig:case2}
\end{figure}

\begin{figure}
\includegraphics[width=2in]{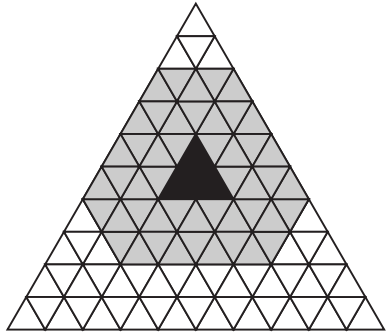}
\caption{Case II: There can be no more than the given point in the greater shaded region, and there is at most one point in the top white triangle. So at least three points lie in the white region at the bottom.}
\label{fig:wide2}
\end{figure}

Case II is that there is a point $p_1$ in the shaded region in Figure \ref{fig:case2}. By symmetry, assume that $p_1$ is in the upper triangle. There are no other points in the wider shaded region of Figure \ref{fig:wide2} or we are done, by the hexagon case of Lemma \ref{trihex}. There is at most one point in the top white triangle, and so there are at least three points in the lower white region. The case of no points in the top triangle and four in the bottom is easy to dispense with by using Lemma \ref{trihex}. So assume that there is one point, $p_2$, in the top triangle.

\begin{figure}
\includegraphics[width=2in]{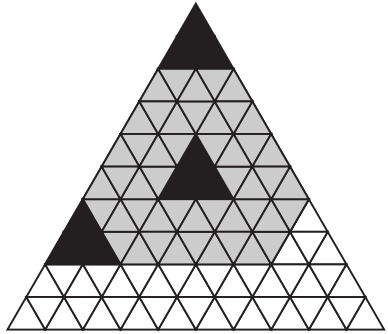}
\caption{Case II:  $p_1$, $p_2$, and $p_3$ are in the three shaded triangles: $p_1$ in the center, $p_2$ at the top, and $p_3$ on the bottom-left.}
\label{fig:case22}
\end{figure}

Now $p_3, p_4,$ and $p_5$ are all in the lower white region of Figure \ref{fig:case2}. They can not all be in the lower two strips, or we are done by the parallelogram lemma with a $(2,10)$-parallelogram, so there is at least one point $p_3$ on the third or fourth row, by symmetry assume it is in the triangle on the left, as in Figure \ref{fig:case22}.

\begin{figure}
\includegraphics[width=2in]{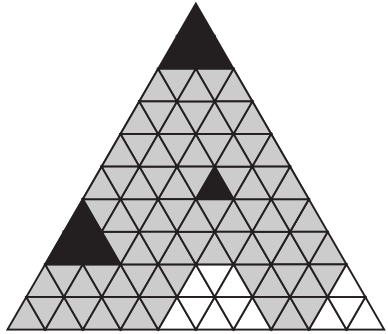}
\caption{Case II: Applying Lemma \ref{parallel} with $(3,8)-$ and $(4,6)-$ parallelograms gives than $p_4$ and $p_5$ are in the white regions at the bottom of $T$. Then $p_3$, $p_4$, and $p_5$ form a small triangle.}
\label{fig:case24}
\end{figure}

By Lemma \ref{parallel} with a $(3,9)$-parallelogram along the left edge of $T$ and containing points $p_2$ and $p_3$, $p_1$ must be in the bottom right corner triangle of side length $L/10$, as in Figure \ref{fig:case24}. Then applying Lemma \ref{parallel} with $(3,8)-$ and $(4,6)$-parallelograms gives that $p_4$ and $p_5$ are in the white region at the bottom of Figure \ref{fig:case24}.

\begin{figure}
\includegraphics[width=2in]{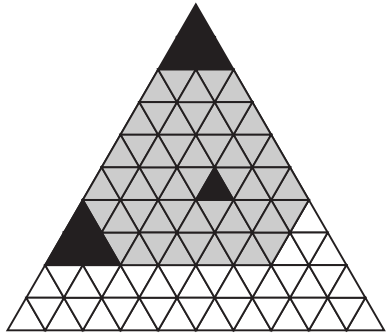}
\caption{Case II: Then by the parallelogram lemma, with a $(3,9)$-parallelogram, point $p_1$ is in the bottom right corner of the center triangle, or we are done.}
\label{fig:case23}
\end{figure}

We claim that then $p_3$, $p_4$, and $p_5$ form a triangle of area $24/100$ or less. Given any particular location of $p_4$ and $p_5$, triangle $p_3p_4p_5$ attains its maximum value (given $p_4$ and $p_5$) with $p_3$ at a corner of its boundary triangle. If $p_3$ is one of the bottom two corners, then the three points are in a $(2,10)$-parallelogram, and if it is the top corner, they are in a $(3,8)$-parallelogram. In either case, the area of the triangle is less than $24/100$ by Lemma \ref{parallel}.

\begin{figure}
\includegraphics[width=2in]{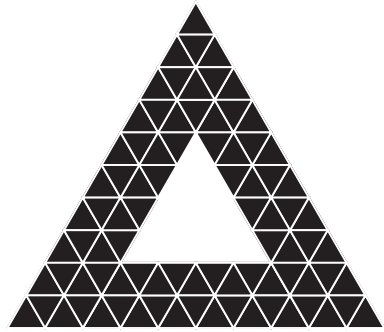}
\caption{Case III: All five points lie in the shaded area.}
\label{fig:case3}
\end{figure}

Then Case III is the only possibility left, that all five points contained in the outer shaded region in Figure \ref{fig:case3}.

\begin{figure}
\includegraphics[width=2in]{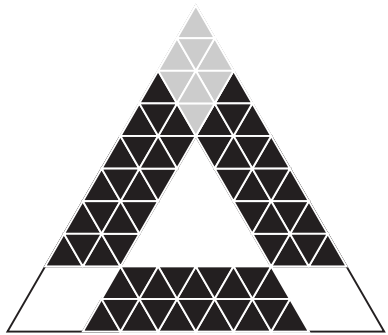}
\caption{Case III: We partition the outer strip into six smaller regions, diamonds and trapezoids. Lemma \ref{trihex} places severe restrictions on where the points may lie, if they are to form no small triangles. In fact without loss of generality we assume that $p_1$ is in the diamond on top, $p_2$ and $p_3$ in the trapezoids on the left and right, and $p_4$ and $p_5$ in the bottom trapezoid.}
\label{fig:case32}
\end{figure}

We break the outer strip into three diamond and trapezoid shaped regions, Figure \ref{fig:case32}.

Any pair of the trapezoids can be covered with a hexagon and two triangles of diameter $4L/10$, so by Lemma \ref{trihex}, they contain a total of at most three points from $S$. It is quick to deduce from this, that the points are distributed $2-1-1$ in the trapezoids, with one point in the diamond opposite the trapezoid with two points, or we are done. By symmetry we assume $p_1$ is in the diamond on top, $p_2$ and $p_3$ are in the trapezoids on the left and right, respectively, and $p_4$ and $p_5$ are in the bottom trapezoid.

\begin{figure}
\includegraphics[width=2in]{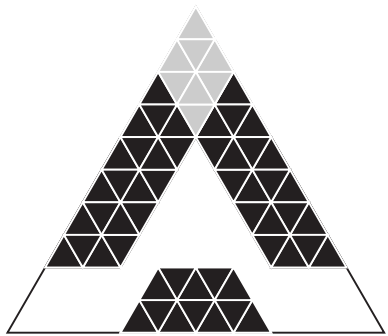}
\caption{Case III: $(3,8)$-parallelograms force $p_4$ and $p_5$ into in an even smaller trapezoid, shaded at bottom, but then we are done by Lemma \ref{trihex}.}
\label{fig:case33}
\end{figure}

By $(3,8)$-parallelograms at the top vertex of $T$ and going down the left and right sides, $p_4$ and $p_5$ are in an even smaller trapezoid, as seen in Figure \ref{fig:case33}, but then we are done by Lemma \ref{trihex}.\\

This completes the proof of the main result, but begs the question:

\begin{question} What is the minimum $s$ of all $\sigma$ such that among every five points in a triangle of unit area, some three of them form a triangle of area less than or equal to $\sigma$? (Note that $s$ is actually a minimum and not just an infinium, by compactness. \cite{Etudes})
\end{question} \label{inf}

\begin{figure}
\includegraphics[width=2.5in]{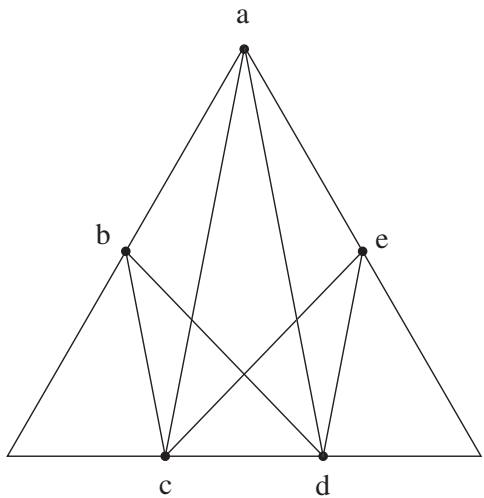}
\caption{Five points in a unit area triangle, with none of the triangles having area less than $1/6$. Point $a$ is at the corner, $b$ and $e$ are on midpoints, and $c$ and $d$ partition the bottom side into three equal segments. Note that triangles $abc$, $ade$, $bcd$, and $cde$ all have area equal to $1/6$.}
\label{fig:s16}
\end{figure}

By Figure \ref{fig:s16}, we have $s \ge 1/6$. This might be best possible.

\begin{conjecture} $s=1/6$.
\end{conjecture}

\section{Many small triangles}

{\bf In this section we return to using ``small'' to refer to triangles of area $\le 1/4$.} There is another way in which Theorem \ref{soif} might be improved, besides varying the $1/4$. Is it possible that given five points in a triangle of unit area, they must form not only one small triangle, but more than one? We conjecture that this is the case; in particular we believe the following holds.

\begin{conjecture} Given five points in a triangle of unit area, then they form at least {\it three} small triangles. \end{conjecture}

\begin{figure}
\includegraphics[width=2in]{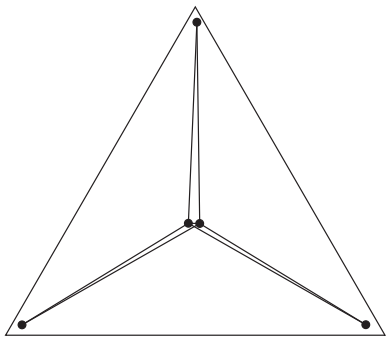}
\caption{There are only three small triangles.}
\label{fig:3small}
\end{figure}

This is not true with ``three'' replaced by ``four'', by the arrangement in Figure \ref{fig:3small}. If the conjecture holds then the following corollary would follow, by averaging: Given $n \ge 5$ points in a triangle of unit area, there they form at least $$ \frac{3}{{5 \choose 2}} {n \choose 3} = \frac{3}{10} {n \choose 3}$$ small triangles. This leads us to the following question.

\begin{question} What is the supremum $A$ of all $\alpha \in [0,1]$ such that given (sufficiently large) $n$ points in a triangle of unit area, they form at least $$ \alpha {n \choose 3}$$ small triangles? \end{question}

\begin{figure}
\includegraphics[width=2in]{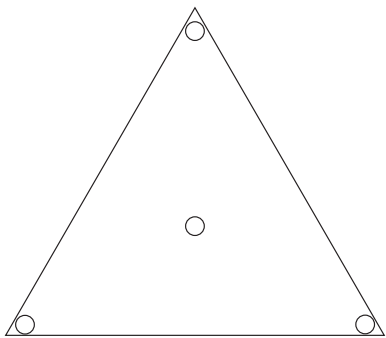}
\caption{If there approximately $n/4$ points in each of the small circles, then there are $4 (n/4)^3 = n^3 / 16 \approx (3/8) {n \choose 3}$ big triangles, hence approximately $(5/8){n \choose 3}$ small triangles. }
\label{fig:58small}
\end{figure}

We observe that $A \le 5/8$, by the construction in Figure \ref{fig:58small}, and conjecture that asymptotically this arrangement gives the smallest number of small triangles.

\begin{conjecture} $A=5/8$. \end{conjecture}

\end{document}